\documentclass[a4paper,11pt]{article}

\usepackage{pifont}
\date{}

\usepackage{makeidx}
\usepackage[centertags]{amsmath}
\usepackage{amsfonts}
\usepackage{amssymb}
\usepackage{amsthm}
\usepackage{accents}
\usepackage{epsfig}
\usepackage{color}

\makeatletter \@addtoreset{equation}{section} \makeatother

\begin{document}

\title{\bf Periodic Solutions of Non-Autonomous Second Order Hamiltonian Systems  \footnote{Supported
by National Natural Science Foundation of China.}}
\author{{ Fengying Li\footnote{Email:lify0308@163.com}} \\
\small \it  The School of Economic and Mathematics, Southwestern
University of Finance and Economics, \\
\small \it Chengdu 611130, China\\
\ Shiqing Zhang\ and \ Xiaoxiao Zhao\\
 \small \it College of
Mathematics, Sichuan University,
 \small\it Chengdu 610064, People's Republic of China}

\maketitle

\begin{quote}

{\bf Abstract}

  We try to generalize a result of M. Willem on forced
periodic oscillations which required the assumption
that the forced potential is periodic on spatial variables. In this
 paper, we only assume its integral on the time variable is
periodic, and so we extend the result to cover the forced pendulum equation.
We apply the direct variational minimizing method and Rabinowtz's saddle point theorem to study the periodic
solution when the integral of the potential on the time variable is
periodic. \\

{\bf Keywords}

 Forced second order Hamiltonian systems, the forced pendulum equation, variational minimizers,  Saddle Point Theorem.\\
 {\bf 2000AMS Subject Classification} 34C15, 34C25.

\end{quote}

\section{Introduction and Main Results}
\ \ \ \ \ In \cite{r3} and [5], M. Willem  and Mawhin studied the following
second order Hamiltonian system
\begin{equation}\label{1}
\ddot{u}(t)=-\nabla F(t,u(t))=-F'(t,u(t))
\end{equation}
where $F: [0,T]\times R^{N}\rightarrow R ,\nabla F(t,u(t))=F'(t,u(t))$
is the gradient of $F(t,u(t))$ with respect to $u$.
We assume $F(t,u(t))$ satisfies the following
assumption:

(A). $F(t,x)$ is measurable in $t$ for each $x\in R^{N}$,
continuously differentiable in $x$ for a.e. $t\in [0,T]$, and there
exist $a\in C(R^{+},R^{+})$ and $b\in L^{1}(0,T;R^{+})$ such that
$$|F(t,x)|\leq a(|x|)b(t),$$
$$|\nabla F(t,x)|\leq a(|x|)b(t)$$ for all
$x\in R^{N}$ and a.e. $t\in [0,T]$.

  M. Willem (\cite{r3}) got the following  theorem :
\vspace{0.4cm}

\textbf{Theorem 1.1}\
(\cite{r3} and [5]) Assume $F$ satisfies
 condition (A) and for the canonical basis $\{e_{i}|1\leq i\leq N\}$
 of $R^{N}$, there exist $T_{i}>0$ such that for $\forall x\in
 R^{N}$ and a.e. $t\in [0,T]$,
\begin{equation}\label{2}
F(t,x+T_{i}e_{i})=F(t,x), \ \ \ \ 1\leq i\leq N
\end{equation}
Then (1.1) has at least one solution which minimizes
$$f(u)=\int_{0}^{T}[\frac{1}{2}|\dot{u}(t)|^{2}-F(t,u(t))]dt$$
on
$H_{T}^{1}=\{u|u, \dot{u}\in L^{2}[0,T], u(t+T)=u(t)\}$.
\vspace{0.4cm}

In order to cover the forced pendulum equation:
\begin{equation}\label{3}
\ddot{u}(t)=-a\sin u+e(t),
\end{equation}

Mawhin-Willem [5] also study the following forced equation:
\begin{equation}\label{1}
\ddot{u}(t)=-\nabla F(t,u(t))-e(t)=-F'(t,u(t))-e(t)
\end{equation}
they got the following Theorem:

\textbf{Theorem 1.2}\
(\cite{r3} and [5]) Assume $F$ satisfies
 the conditions of Theorem 1.1,and $e(t)\in L^{1}(0,T;R^N)$ verifying
 $$\int_{0}^{T}e(t)dt=0,$$
 then (1.4) has at least one solution which minimizes on $H_{T}^{1}$ the following functional:
$$f(u)=\int_{0}^{T}[\frac{1}{2}|\dot{u}(t)|^{2}-F(t,u(t))-e(t)u(t)]dt$$
   We notice that the potential $F(t,x)=-(a\cos x+e(t)x)$ does not satisfy (\ref{2}). But
if $\int_{0}^{T}e(t)dt=0$, then $F(t,x)=-(a\cos x+e(t)x)$ does satisfy
\begin{equation}\label{4}
\int_{0}^{T}F(t,x+2\pi)dt=\int_{0}^{T}F(t,x)dt.
\end{equation}
So instead of (\ref{2}) we only assume the weaker integral condition:
\begin{equation}\label{5}
\int_{0}^{T}F(t,x+T_{i}e_{i})dt=\int_{0}^{T}F(t,x)dt \ \ \ \
i=1,2,...,N
\end{equation}
We obtain the following results:

 \vspace{0.4cm}\textbf{Theorem 1.3}\ Assume $F: R\times R^{N}\rightarrow
 R$ satisfies condition (A) and\\
 (F1). $F(t+T,x)=F(t,x)$, $\forall (t,x)\in R\times R^{N}$,\\
 (F2). $F$ satisfies (\ref{5}).\\
 (F3). There exist $0<C_{1}<\frac{1}{2}(\frac{2\pi}{T})^2$, $C_2>0$ such that
 $$|F(t,x)|\leq C_1|x|^2+C_2$$

  Then (1.1) has at least one $T$-periodic
 solution.

\vspace{0.4cm}

\textbf{Corollary 1.1} (J. Mawhin, M. Willem
\cite{r5})\ For the pendulum equation (\ref{3}), the potential
$F(t,x)=a\cos x+e(t)x$ satisfies all conditions in Theorem 1.3
provided $e(t+T)=e(t)$ and $\int_{0}^{T}e(t)dt=0$. In this case,
(\ref{3}) has at least one $T$-periodic solution.

\vspace{0.4cm}

\textbf{Theorem 1.4}\ Suppose $F: R\times
R^N\rightarrow R$ satisfies conditions (A), (F1), (F2) and\\
 (F4).
There are $\mu_1<2$, $\mu_2\in R$ such that
$$F'(t,x)\cdot x\leq \mu_1F(t,x)+\mu_2,$$\\
(F5). There is $\delta>0$ such that for $ t\in R$, $F(t,x)>\delta$, as
$|x|\rightarrow +\infty$,\\
(F6). $F(t,x)\leq b|x|^2$.\\

 Then if $T<\sqrt{\frac{2}{b}}\pi$, (1.1) has a $T-$periodic solution;
furthermore, if $\forall x\in R^N$, $\int_0^T F(t,x)dt\geq 0$, then
(\ref{1}) has a non-constant $T-$periodic solution.

\section{Some Important Lemmas}

\ \ \ \ \ \ \ \ {\bf Lemma 2.1}\ (Eberlin-Smulian\cite{ES}) A
Banach space X is reflexive if and only if any bounded sequence in X
has a weakly convergent subsequence.
\vspace{0.3cm}

 {\bf Lemma 2.2}\
(\cite{e1},\cite{r1},\cite{e15})\ Let $q\in W^{1,2}(R/TZ,R^{n})$ and
$\int^T_0q(t)dt=0$, then

(i). We have Poincare-Wirtinger's inequality
$$\int^T_0|\dot{q}(t)|^2dt\geq (\frac{2\pi}{T})^2\int_0^T|q(t)|^2dt$$

(ii). We have Sobolev's inequality
$$\max_{0\leq t\leq T}|q(t)|=\|q\|_{\infty}\leq
\sqrt{\frac{T}{12}}(\int^T_0|\dot{q}(t)|^2dt)^{1/2}$$
\vspace{0.2cm}

We define the equivalent norm in $H^{1}_{T}=H^{1}=W^{1,2}(R/TZ,R^n):$
$$\|q\|_{H^1}=(\int_0^T|\dot{q}|^2)^{1/2}+|\int_0^T q(t)dt|$$

{\bf Lemma 2.3}(\cite{r1})\ Let $X$ be a reflexive Banach space,
$M\subset X$ a weakly closed subset, and $f: M\rightarrow R\cup
\{+\infty\}$ weakly lower semi-continuous. If the minimizing
sequence for $f$ on $M$ is bounded, then $f$ attains its infimum on
$M$.
\vspace{0.4cm}

{\bf Definition 2.1}(\cite{e4})\ Suppose $X$ is a Banach space and $f\in C^{1}(X,R)$ and
$\{q_{n}\}\subset X$ satisfies
$$f(q_{n})\rightarrow C,\ \ \ \ (1+\|q_{n}\|)f'(q_{n})\rightarrow 0.$$
Then we say $\{q_n\}$ satisfies the $(CPS)_C$ condition.
\vspace{0.4cm}

 {\bf Lemma 2.4}(Rabinowitz's Saddle Point Theorem\cite{e12},
Mawhin-Willem\cite{r1})\ Let $X$ be a Banach space with $f\in
C^1(X,R)$. Let $X=X_1\oplus X_2$ with
$$dim X_1<+\infty$$
and
$$\sup_{S^1_R} f<\inf_{X_2} f,$$
where $S^1_R=\{u\in X_1||u|=R\}$.

Let $B_R^1=\{u\in X_1, |u|\leq R\}$, $M=\{g\in C(B^1_R,X)|g(s)=s,
s\in S^1_R\}$
$$C=\inf_{g\in M}\max_{s\in B^1_R}f(g(s)).$$
Then $C>\inf_{X_2}f$, and if $f$ satisfies $(CPS)_{C}$ condition, then $C$ is
a critical value of $f$.
\vspace{0.4cm}

\section{The Proof of Theorem 1.3}

\vspace{0.4cm}

\textbf{Lemma 3.1} (Morrey \cite{r2}, M-W \cite{r3})\
Let $L: [0,T]\times R^{N}\times R^{N}\rightarrow R$,
$(t,x,y)\rightarrow L(t,x,y)$ be measurable in $t$ for each
$(x,y)\in R^{N}\times R^{N}$ and continuously differentiable in
$(x,y)$ for a.e. $t\in [0,T]$. Suppose there exists $a\in
C(R^{+},R^{+})$, $b\in L^{1}(0,T; R^{+})$ and $c\in L^{q}(0,T;
R^{+})$, $1<q<\infty$, such that for a.e. $t\in [0,T]$ and every
$(x,y)\in R^{N}\times R^{N}$ one has
$$|L(t,x,y)|\leq a(|x|)(b(t)+|y|^{q}),$$
$$|D_{x}L(t,x,y)|\leq a(|x|)(b(t)+|y|^{p}),$$
$$|D_{y}L(t,x,y)|\leq
a(|x|)(c(t)+|y|^{p-1}).$$
 where $\frac{1}{p}+\frac{1}{q}=1$, then the
functional
$$\varphi(u)=\int_{0}^{T}L(t,u(t),\dot{u}(t))dt$$
is continuously differentiable on the Sobolev space
\begin{equation}\label{6}
W^{1,p}=\{u\in L^{p}(0,T), \dot{u}\in L^{p}(0,T)\}
\end{equation}
and
\begin{equation}\label{7}
<\varphi'(u),v>=\int_{0}^{T}[<D_{x}L(t,u,\dot{u}),v>+D_{y}L(t,u,\dot{u})\cdot\dot{v}]dt
\end{equation}

From Lemma 3.1 and the assumptions (A), we know that the variational
functional
\begin{equation}\label{8}
f(u)=\int_{0}^{T}[\frac{1}{2}|\dot{u}|^{2}-F(t,u(t))]dt
\end{equation}
is $C^{1}$ on $W_{T}^{1,2}=H^{1}_{T}$, and the critical point is
just the periodic solution for the system (1.1).

Furthermore, if (F1) and (F2) are satisfied, we will prove the
functional $f(u)$ attains its infimum on $H_{T}^{1}$; in fact,
\begin{equation}\label{10}
H_{T}^{1}=X\oplus R^{N},
\end{equation}
where
\begin{equation}\label{11}
 X=\{x\in H^{1}_{T}: \bar{x}\triangleq\frac{1}{T}\int_{0}^{T}x(t)dt=0\}
\end{equation}
 and $\forall u\in H_{T}^{1}$, we have $\widetilde{u}\in X$ and
$\overline{u}\in R^{N}$, such that $u=\widetilde{u}+\overline{u}$.
\vspace{0.2cm}

By Poincare-Wirtinger's inequality,
\begin{equation}\label{12}
  \begin{aligned}
f(\tilde{u})&\geq\frac{1}{2}\int_{0}^{T}|\dot{\tilde{u}}|^{2}dt-C_{1}\int_{0}^{T}|\tilde{u}|^{2}dt-C_2T\\
&\geq[\frac{1}{2}-C_1(\frac{2\pi}{T})^{-2}]\int_{0}^{T}|\dot{\widetilde{u}}|^{2}dt-C_2T;
 \end{aligned}
\end{equation}
hence, $f$ is coercive on $X$.

Let $\{u_{k}\}$ be a minimizing sequence for $f(u)$ on $H_{T}^{1}$,
$u_{k}=\widetilde{u}_{k}+\overline{u}_{k}$, where
$\widetilde{u}_{k}\in X$, $\overline{u}_{k}\in R^{N}$, then by
(\ref{12}) we have
\begin{equation}\label{14}
\|\widetilde{u}_{k}\|_{H^{1}_{T}}\leq C.
\end{equation}
By condition (F2), we have
\begin{equation}\label{15}
f(u+T_{i}e_{i})=f(u), \ \ \ \ \forall u\in H_{T}^{1}, \ \ \ \ 1\leq
i\leq N.
\end{equation}
So if $\{u_{k}\}$ is a minimizing sequence for $f$, then
$$(\widetilde{u}_{k}\cdot e_{1}+\overline{u}_{k}\cdot e_{1}+k_{1}T_{1},...,\widetilde{u}_{k}\cdot e_{N}+\overline{u}_{k}\cdot e_{N}+k_{N}T_{N})$$
is also a minimizing sequence of $f(u)$, and so we can assume
\begin{equation}\label{16}
0\leq\overline{u}_{k}\cdot e_{i}\leq T_{i}, \ \ \ \ 0\leq i\leq N.
\end{equation}
By (\ref{14}) and (\ref{16}), we know $\{u_{k}\}$ is a bounded
minimizing sequence in $H_{T}^{1}$, and it has a weakly convergent
subsequence; furthermore, $f$ is weakly lower semi-continuous since
$f$ is the sum of a convex continuous function and a weakly
continuous function. We can conclude that $f$ attains its infimum on
$H_{T}^{1}$. The
corresponding minimizer is a periodic solution of (\ref{1}).\\

\section{The Proof of Theorem 1.4}

{\bf Lemma 4.1}:\ If conditions (A), (F1), (F2) and (F4) in Theorem
1.4 hold, then $f(q)$ satisfies the $(CPS)_C$ condition on $H^1$.

{\bf Proof}: For any $C$, let $\{u_n\}\subset H^1$ satisfy
\begin{equation}\label{17'}
f(u_n)\rightarrow C, \ \ \ \ (1+\|u_n\|)f'(u_n)\rightarrow 0.
\end{equation}
We claim $\|\dot{u}_n\|_{L^2}$ is bounded; in fact, by
$f(u_n)\rightarrow C$, we have
\begin{equation}\label{18'}
\frac{1}{2}\|\dot{u}_n\|^2_{L^2}-\int_0^T F(t,u_n)dt\rightarrow C.
\end{equation}
By (F4) we have
\begin{equation}\label{19'}
\begin{aligned}
 <f'(u_n),u_n>&=&\|\dot{u}_n\|^2_{L^2}-\int_0^T
 (<F'(t,u_n),u_n>)dt\\
 &\geq&\|\dot{u}_n\|^2_{L^2}-\int_0^T [\mu_2+\mu_1F(t,u_n)]dt.
 \end{aligned}
\end{equation}
By (\ref{18'}) and (\ref{19'}), we see that

\begin{equation}\label{20'}
0\leftarrow<f'(u_n),u_n>\geq a\|\dot{u}_n\|^2_{L^2}+C_1+\delta,
n\rightarrow +\infty
\end{equation}
where $C_1=C\mu_1-T\mu_2+\delta, \delta>0, a=1-\frac{\mu_1}{2}>0.$\\
We have shown that $\|\dot{u}_n\|_{L^2}$ is bounded.

By condition (F2), we have
\begin{equation}\label{15''}
f(u+T_{i}e_{i})=f(u), \ \ \ \ \forall u\in H_{T}^{1}, \ \ \ \ 1\leq
i\leq N.
\end{equation}
Hence, if $\{u_{k}\}$ is a $(CPS)_C$ sequence for $f$, then
$$(\widetilde{u}_{k}\cdot e_{1}+\overline{u}_{k}\cdot e_{1}+k_{1}T_{1},...,\widetilde{u}_{k}\cdot e_{N}+\overline{u}_{k}\cdot e_{N}+k_{N}T_{N})$$
is also a $(CPS)_C$ sequence of $f(u)$, so we can assume
\begin{equation}\label{16''}
0\leq\overline{u}_{k}\cdot e_{i}\leq T_{i}, \ \ \ \ 0\leq i\leq N.
\end{equation}
By (\ref{16''}), we know $|\bar{u}_{k}|$ is bounded, and so $\|u_n\|=\|\dot{u}_n\|_{L^2}+|\int_0^T u_n(t)dt|$ is bounded.

The rest of he lemma can be completed in a now standard fashion.

We finish the proof of \textbf{Theorem 1.4.} In Rabinowitz's Saddle Point
Theorem, we take
$$X_1=R^N, X_2=\{u\in W^{1,2}(R/TZ,R^N),\int_0^T u dt=0\}.$$

For $u\in X_2$, we may use the Poincare-Wirtinger inequality, and so by Lemma
2.2 and (F6), we have
\begin{eqnarray*}\label{22}
f(u)&\geq&\frac{1}{2}\int_0^T|\dot{u}|^2 dt-b\int_0^T|u|^2 dt\\
&\geq&[\frac{1}{2}-b(2\pi)^{-2}T^2]\int_0^T|\dot{u}|^2dt\\
&\geq& 0.
\end{eqnarray*}
On the other hand, if $u\in R^N$, then by $(F5)$ we have
$$f(u)=-\int_0^T F(t,u)dt\leq -\delta, |u|=R\rightarrow
+\infty.$$

The proof of Theorem 1.4 is concluded by calling upon Rabinowitz's Saddle Point Theorem.
In fact,there is a critical point $\bar{u}$ such that $f(\bar{u})=C>\inf_{X_2} f(u)\geq 0$,
which is nonconstant since otherwise $f(\bar{u})=-\int_0^TF(\bar u,t)dt\leq 0,$ which is a contradiction.\\

The authors would like to thank the referees for their valuable suggestions.

\end{document}